\documentclass[a4paper, 12pt]{article}
\usepackage{amsmath}
\usepackage{mathrsfs}
\usepackage{latexsym,epsfig,amsmath,bm}
\usepackage{amssymb}
\usepackage{ebezier}
\usepackage{tikz}
\usepackage{epic,eepic,color,xcolor,caption}
\usetikzlibrary{decorations.pathreplacing}
\newcounter{prostredi}
\def\theprostredi{\arabic{prostredi}}

\makeatletter
\def\@begintheorem#1#2{\trivlist
   \item[\hskip \labelsep{\bfseries #1\ #2.}]} 
\def\@opargbegintheorem#1#2#3{\trivlist
      \item[\hskip \labelsep{\bfseries #1\ #2.\ (#3.)}]} 
\def\@endtheorem{\endtrivlist}
\makeatother



\newenvironment{claim}{\par\bigskip\noindent%
\refstepcounter{prostredi}{\bf Claim \theprostredi.}\quad\bgroup\sl
}
{\egroup\par\bigskip\endtrivlist}%

\newcommand{\qed}{{$\qquad\square$\bigbreak}}

\newtheorem{theorem}{Theorem}
\newtheorem{lemma}[theorem]{Lemma}

\newtheorem{corollary}[theorem]{Corollary}

\tikzstyle{fleche}=[->,>=stealth',thick,rounded corners=4pt]

\title{On traceable iterated line graph and hamiltonian path index}

\author{\small{Zhaohong Niu\thanks{School of Mathematical Sciences, Shanxi University, Taiyuan, Shanxi, 030006, P.R. China. E-mail address: zhniu@sxu.edu.cn. Research is supported by the Natural Science Funds for Young Scientists of China (No. 11501341 and 11701349).}\,\,,
Liming Xiong\thanks{School of Mathematics and Statistics, Beijing Institute of
Technology, Beijing, 100081, P.R. China. Research is supported by the Natural Science Funds of China (No. 11471037).}\,\,,Weihua Yang\thanks{Department of Mathematics, Taiyuan University of Technology, Taiyuan, Shanxi, 030024, P.R. China. Research is supported by the Natural Science Funds of China (No. 11671296).}}}

\date{}
\begin{document}
\maketitle

\begin{abstract}
Xiong and Liu [L. Xiong and Z. Liu, Hamiltonian iterated line graphs, Discrete Math. 256 (2002) 407-422] gave a characterization of the graphs $G$ for which the $n$-th iterated line graph $L^n(G)$ is hamiltonian, for $n\ge2$. In this paper, we study the existence of a hamiltonian path in $L^n(G)$, and give a characterization of $G$ for which $L^n(G)$ has a hamiltonian path. As applications, we use this characterization to give several upper bounds on the hamiltonian path index of a graph.\\


\end{abstract}

\normalsize

\section{Introduction}\label{sec:intro}

\hspace*{\parindent}The graphs considered in this paper are finite, undirected and loopless. Undefined notation and terminology will follow
\cite{Bondy}.

The \textit{line graph} $L(G)$ of a graph $G$ has $E(G)$ as its vertex set and two vertices are adjacent in $L(G)$ if and only if they are adjacent as edges in $G$. A trail $T$ of $G$ is \textit{dominated} if each edge of $G$ is incident with at least one vertex of $T$.

\begin{theorem}(Harary and Nash-Williams, \cite{HN})\label{th:HN} \itshape Let $G$ be a connected graph with at least three edges. Then $L(G)$ is hamiltonian if and only if $G$ has a dominating closed trail.\end{theorem}

Theorem \ref{th:HN} characterized those graphs $G$ (expecially, non-hamiltonian graphs) for which $L(G)$ is hamiltonian. It is natural to believe that for most graphs, after applying the line graph operation iteratively a finite  number of times, the resulting will be hamiltonian. Chartrand \cite{C} showed this is indeed the case. He considered the $n$-iterated line graph $L^n(G)$ of $G$, and introduced the \textit{hamiltonian index} of a graph, denoted by $h(G)$, i.e., the minimum number $n$ such that $L^n(G)$ is hamiltonian. Here the $n$-\textit{iterated line graph} of a graph $G$ is defined to be $L(L^{n-1}(G))$, where $L^1(G)$ denotes the line graph $L(G)$ of $G$, and $L^{n-1}(G)$ is assumed to have a nonempty edge set. Chartrand showed that for any graph $G$ other than a path, the hamiltonian index of $G$ exists. Since then, many results of determing the exact values  or upper bounds of $h(G)$ for special graphs are proved, such as Chartrand and Wall \cite{CW} for trees (other than paths) and connected graphs $G$ with $\delta(G)\ge3$, Kapoor and Stewart \cite{KS} for $K_{2,n}\,(n\ge3)$, Xiong and Liu \cite{X} for connected graphs using the split blocks and Catlin's reduction method, Han et al. \cite{Han} for 2-connected graph with $\kappa(G)\ge\alpha(G)-t$ (where $t$ is a nonnegative integer). Ryj\'a\v{c}ek et al. \cite{Ry} showed that the problem to decide whether the hamiltonian index of a
given graph is less than or equal to a given constant is NP-complete.

Let $G$ be a graph and $H$ a subgraph of $G$, then $V(H)$ and $E(H)$ denote the sets of vertices and edges of $H$, respectively. Define $V_i(H)=\{v\in V(H):d_H(v)=i\}$ and $W(H)=V(H)\backslash V_2(H)$. A \textit{branch} in $G$ is a nontrivial path with ends in $W(G)$ and with internal vertices, if any, of degree 2. We denote by $B(G)$ the set of branches of $G$. Define $B_1(G)=\{b\in B(G):V(b)\cap V_1(G)\neq\emptyset\}$. The distance between two subgraphs $H_1$ and $H_2$ of $G$, denoted by $d_G(H_1,H_2)$, is $\min\{d_G(v_1,v_2): v_1\in V(H_1)\text{ and }v_2\in V(H_2)\}$.

For a positive integer $k$, Xiong and Liu \cite{X} used $EU_k(G)$ to denote the set of subgraphs $H$ of a graph $G$ that satisfy the following conditions, and then investigated the existence of a characterization of those graphs $G$ for which $L^n(G)$ is hamiltonian, which was mentioned in \cite{CM}, and proved Theorem \ref{th:X:main}.

\begin{itemize}
\item[(I)] $d_H(x)\equiv0\,(\!\!\!\!\mod2)$ for every $x\in V(H)$;
\item[(II)] $V_0(H)\subseteq \bigcup_{i=3}^{\Delta(G)}V_i(G)\subseteq V(H)$;
\item[(III)] $d_G(H_1,H-H_1)\le k-1$ for every subgraph $H_1$ of $H$;
\item[(IV)] $|E(b)|\le k+1$ for every branch $b\in B(G)$ with $E(b)\cap E(H)=\emptyset$;
\item[(V)] $|E(b)|\le k$ for every branch $b\in B_1(G)$.
\end{itemize}

\begin{theorem}(Xiong and Liu, \cite{X})\label{th:X:main} \itshape Let $G$ be a connected graph with at least three edges and $n\ge2$. Then $L^n(G)$ is hamiltonian if and only if $EU_n(G)\neq\emptyset$.\end{theorem}

With the help of Theorem \ref{th:X:main}, one can deduce $h(G)\le k$ if it is convenient to check that $EU_k(G)$ is nonempty. See \cite{HX,XB,XW} for more upper bounds of $h(G)$.

The further research of hamiltonianity has been producing several related branches. As a form of weakening, the existence of hamiltonian paths of a graph $G$ or a line graph $L(G)$ is also got a lot of attention. A graph $G$ is \textit{traceable} if it has a hamiltonian path.  Mom\`ege \cite{Mo} showed that a connceted graph $G$ with  $\sigma_2(G)\ge2n/3$ and $K_{1,4}$-free is traceable. He and Yang \cite {He} proved that there exist at least $\max\{1,\lfloor\frac{\,1\,}{\,8\,}\delta(G)\rfloor-1\}$ edge-disjoint hamiltonian paths between any two vertices in a hamiltonian-connected line graph $L(G)$. Xiong and Zong \cite{XZ} discussed the traceability of line graphs, and obtained a similar result as Theorem \ref{th:HN}. Note that the following result is also an immediately corollary of Proposition 2.2 in \cite{Lai}, in which Lai et al. discussed the  hamiltonian connectedness of $L(G)$. 

\begin{theorem}(Xiong and Zong, \cite{XZ}, Lai, Shao and Zhan, \cite{Lai})\label{th:zong} \itshape Let $G$ be a connected graph with at least three edges. Then the line graph $L(G)$ is traceable if and only if $G$ has a dominating trail.\end{theorem}

In this paper, inspired by Theorems \ref{th:HN}, \ref{th:X:main} and \ref{th:zong}, and the closed relationship between hamiltonian graphs and traceable graphs, we study the existence of a hamiltonian path in the $n$-iterated line graph $L^n(G)$, i.e., the traceability of $L^n(G)$. 

For a graph $G$, we use $O(G)$ to denote the set of odd degree vertices of $G$. The condition (I) guarantees that each vertex in $V(H)$ has even degree, i.e., $O(H)=\emptyset$. When concerning about the traceability of $L^n(G)$, this is too strong. That is, we can allow that $H$ has at most two odd degree vertices. Moreover, by the fact that $H$ may contain a branch $b\in B_1(G)$ or a part of $b$, we just need the condition (V) holds for every branch $b\in B_1(G)$ with $E(b)\cap E(H)=\emptyset$.

Let $EUP_k(G)$ denote the set of subgraphs $H$ of a graph $G$ that satisfy the following conditions (I)$'$ and (V)$'$, and the conditions (II)-(IV) in the definition of $EU_k(G)$.

\begin{itemize}
\item[(I)$'$] $|O(H)|\le2$;
\item[(V)$'$] $|E(b)|\le k$ for every branch $b\in B_1(G)$ with $E(b)\cap E(H)=\emptyset$.
\end{itemize}

Obviously, $EU_k(G)\subseteq EUP_k(G)$. Now we are ready to present our main result, the proof of which will be postponed to Section \ref{sec:proof}.

\begin{theorem}\label{th:the main} \itshape Let $G$ be a connected graph with at least three edges and $n\ge2$. Then $L^n(G)$ is traceable if and only if $EUP_n(G)\neq\emptyset$.\end{theorem}

Theorem \ref{th:the main} doesn't hold for $n=1$. For example, Fig. 1 shows a graph $G$ with $EUP_1(G)=\emptyset$ while $L(G)$ is traceable. By the definition of $EUP_1(G)$, any subgraph $H$ in $EUP_1(G)$ should be connected, $V_3(G)=\{v_3,v_6,v_9,v_{12}\}\subseteq V(H)$, and the 4 branches $v_1v_2v_3$, $v_3v_4v_5v_6$, $v_6v_7v_8v_9$ and $v_9v_{10}v_{11}$ belong to $H$. Then $|O(H)|\ge4$, a contradiction to (I)$'$. Thus, $EUP_1(G)=\emptyset$, but $L(G)$ is traceable by the fact that $v_1v_2\cdots\,\! v_{11}$ is a dominating trail of $G$ and by Theorem \ref{th:zong}.

\vspace{.2in}
\begin{center}
\begin{tikzpicture}[line width=.6pt]
\draw (8,0) arc (-40:-140:5.216);
\filldraw (0,0) circle [radius=1.8pt]
            (0.64,-.64) circle [radius=1.8pt]
            (1.39,-1.165) circle [radius=1.8pt]
            (2.205,-1.546) circle [radius=1.8pt]
            (3.094,-1.785) circle [radius=1.8pt]
            (4,-1.862) circle [radius=1.8pt]
            (4,0) circle [radius=1.8pt]
            (4.906,-1.785) circle [radius=1.8pt]
            (5.795,-1.546) circle [radius=1.8pt]
            (6.61,-1.165) circle [radius=1.8pt]
            (7.36,-.64) circle [radius=1.8pt]
            (8,0) circle [radius=1.8pt];
\draw (0,0) node[below]{$v_1$};
\draw (0.64,-.64) node[below]{$v_2$};
\draw (1.39,-1.165) node[below]{$v_3$};
\draw (2.205,-1.546) node[below]{$v_4$};
\draw (3.094,-1.785) node[below]{$v_5$};
\draw (4.05,-1.862) node[below]{$v_6$};
\draw (4.936,-1.785) node[below]{$v_7$};
\draw (5.85,-1.546) node[below]{$v_8$};
\draw (6.71,-1.165) node[below]{$v_9$};
\draw (7.56,-.64) node[below]{$v_{10}$};
\draw (8.2,0) node[below]{$v_{11}$};
\draw (4,0) node[above]{$v_{12}$};
\draw (4,0) -- (4,-1.862);
\draw (4,0) -- (1.39,-1.165);
\draw (4,0) -- (6.61,-1.165);
\end{tikzpicture}

\vspace{.1in}
\small Fig. 1. A graph $G$ with $EUP_1(G)=\emptyset$ while $L(G)$ is traceable.
\end{center}

Moreover, we also examine the \textit{hamiltonian path index} of a graph, denoted by $h_p(G)$, i.e., the minimum number $n$ such that $L^n(G)$ is traceable. Regard $h_p(G)=0$ if $G$ is traceable.

In Section \ref{sec:aux}, we will present more terminology and notation, and some auxiliary results. Theorem \ref{th:the main} will be proved in Section \ref{sec:proof}. As applications, we will give some upper bounds of the hamiltonian path index $h_p(G)$ in the last section.

\section{Auxiliary results}\label{sec:aux}

\hspace*{\parindent}In this section, we present several auxiliary results, which will be used in the proof of Theorem \ref{th:the main}. 

The multi-graph of order 2 with two edges will be called 2-\textit{cycle}. Let $G$ be a graph and $H$ a subgraph of $G$, then $\bar{E}(H)$ denotes the set of all edges of $G$ that are incident with vertices of $H$. If $u\in V(H)$, then $E_H(u)$ denotes the set of all edges of $H$ that are incident with $u$, and $d_H(u)=|E_H(u)|$. A graph is called a \textit{circuit} if it is connected and every vertex has an even degree. Regard $K_1$ as a circuit.

For any subgraph $C$ of $L(G)$, by $S(G,C)$ we denote the collection of circuits $H$ of $G$, such that $L(G[\bar{E}(H)])$ contains $C$, and $C$ contains all elements of $E(H)$. Here and throughout, $G[S]$ denotes the subgraph of $G$ induced by $S$, where $S\subseteq V(G)$ or $S\subseteq E(G)$.

\begin{lemma}(Xiong and Liu, \cite{X})\label{lem:x:8} \itshape Each of the following holds.
\begin{itemize}
\item[(1)] If $C$ is a cycle of $L(G)$ with $|E(C)|\ge3$, then $S(G,C)$ is nonempty.
\item[(2)] If $G$ has a circuit $H$ such that $\bar{E}(H)$ has at least three edges, then $L(G)$ has a cycle $C$ with $V(C)=\bar{E}(H)$.
\end{itemize}
\end{lemma}

\begin{lemma}(Beineke, \cite{B})\label{lem:claw-free} \itshape $K_{1,3}$ is not an induced subgraph of the line graph of any graph.
\end{lemma}

The following lemma indicates the relationship between a branch of $G$ and the corresponding branch of $L(G)$.

\begin{lemma}(Xiong and Liu, \cite{X})\label{lem:x:10} \itshape Let $b=u_1u_2\cdots u_s\,(s\ge3)$ be a path of $G$ and $e_i=u_iu_{i+1}$. Then $b\in B(G)$ if and only if $b'=e_1e_2\cdots e_{s-1}\in B(L(G))$.
\end{lemma}

\begin{lemma}\label{lem:x:11} \itshape Let $H$ be a subgraph of $G$ in $EUP_k(G)$ with a minimum number of components. Then there exists no multiple edges in $\bar{E}(H_1)\cap \bar{E}(H_2)$ for any two components $H_1$ and $H_2$ of $H$.
\end{lemma}

A similar result as Lemma \ref{lem:x:11} was proved by Xiong and Liu \cite{X} for $H\in EU_k(G)$. Then arguing similarly, one can obtain Lemma \ref{lem:x:11}. Hence, we omit the details here. An \textit{eulerian subgraph} of $G$ is a circuit which contains at least one cycle of length at least 3.

\begin{lemma}(Xiong and Liu, \cite{X})\label{lem:x:12} \itshape Let $G$ be a connected graph and $C$ be an eulerian subgraph of the line graph $L(G)$. Then there exists a subgraph $H$ of $G$ with
\begin{itemize}
\item[(1)] $d_H(x)\equiv0\,(\!\!\!\mod2)$ for every $x\in V(H)$;
\item[(2)] $d_G(x)\ge3$ for every $x\in V(H)$ with $d_H(x)=0$;
\item[(3)] for any two components $H^0, H^{00}$ of $H$, there exists a sequence of components $H^0=H_1,H_2,\dots,H_s=H^{00}$ of $H$ such that $d_G(H_i,H_{i+1})\le1$ for $i\in\{1,2,\dots,s-1\}$;
\item[(4)] $L(G[\bar{E}(H)])$ contains $C$, and $C$ contains all elements of $E(H)$.
\end{itemize}
\end{lemma}

\begin{lemma}\label{lem:x:8'} \itshape Each of the following holds.
\begin{itemize}
\item[(1)] If $P$ is a non-trivial path of $L(G)$, then $G$ has a trail $T'$, such that $L(G[\bar{E}(T')])$ contains $P$, and $P$ contains all elements of $E(T')$.
\item[(2)] If $G$ has a connected subgraph $H$ such that $|O(H)|=2$, then $L(G)$ has a path $P$ with $V(P)=\bar{E}(H)$.
\end{itemize}
\end{lemma}

\medbreak\noindent \textbf{Proof.} (1) The proof just need a slight modification of the proof of Theorem \ref{th:HN} in \cite{HN}. So we omit the details here.

(2) Suppose $O(H)=\{u,v\}$. If $|\bar{E}(H)|=1$, then $G\cong H\cong K_2$, (2) holds trivially. So we may assume that $\bar{E}(H)$ has at least two edges. Let $e^*=uv$ be a new edge, which doesn't belong to $E(G)$. Note that if $uv\in E(G)$, then $e^*$ and $uv\,(\in E(G))$ are multiple edges in $G+e^*$. Hence, $H+e^*$ is a circuit of $G+e^*$ such that $\bar{E}(H+e^*)$ has at least three edges. By Lemma \ref{lem:x:8} (2), $L(G+e^*)$ has a cycle $C$ with $V(C)=\bar{E}(H+e^*)$. Let $P=C-v_{e^*}$, where $v_{e^*}$ is the vertex in $L(G+e^*)$ corresponding to the edge $e^*$ in $G$. Note that $\bar{E}(H)\cup\{e^*\}=\bar{E}(H+e^*)$. $P$ is a path of $L(G)$ with $V(P)=\bar{E}(H)$.

So Lemma \ref{lem:x:8'} holds.\qed

In \cite{X}, Theorem \ref{th:X:main} was derived from the following two results by induction. The first one indicates a close relationship between $EU_k(L(G))$ and $EU_{k+1}(G)$, and the second one characterizes graphs with 2-iterated line graphs are hamiltonian. 


\begin{theorem}(Xiong and Liu, \cite{X})\label{th:x:induction step} \itshape Let $G$ be a connected graph and $k\ge1$ be an integer. Then $EU_k(L(G))\neq\emptyset$ if and only if $EU_{k+1}(G)\neq\emptyset$.\end{theorem}

\begin{theorem}(Xiong and Liu, \cite{X})\label{th:x:basis step} \itshape Let $G$ be a connected graph with at least three edges. Then $L^2(G)$ is hamiltonian if and only if $EU_2(G)\neq\emptyset$.\end{theorem}

\section{Proof of Theorem \ref{th:the main}}\label{sec:proof}

\hspace*{\parindent}In this section, we will prove Theorem \ref{th:the main}, which is a direct consequence of the following two theorems. The \textit{symmetric difference} of two non-empty sets $A$ and $B$, denoted by $A\Delta B$, is the set $(A\cup B)\backslash (A\cap B)$.

\begin{theorem}\label{th:induction step} \itshape Let $G$ be a connected graph and $k\ge1$ be an integer. Then $EUP_k(L(G))\neq\emptyset$ if and only if $EUP_{k+1}(G)\neq\emptyset$.\end{theorem}

\medbreak\noindent \textbf{Proof.} \textit{Sufficiency}. Supposing that $EUP_{k+1}(G)\neq\emptyset$. Note that if $EU_{k+1}(G)\neq\emptyset$, then by Theorem \ref{th:x:induction step}, $EU_k(L(G))\neq\emptyset$, and hence, $EUP_k(L(G))\neq\emptyset$ by the fact that $EU_k(G)\subseteq EUP_k(G)$. So we may assume that $EU_{k+1}(G)=\emptyset$, which implies that each subgraph of $G$ in $EUP_{k+1}(G)$ contains exactly two odd vertices.

Now let $H\in EUP_{k+1}(G)$ with a minimum number of components denoted by $C_1,C_2,\dots,C_t$. If $H$ is connected, then $H\cong C_1$. Without loss of generality, we let $|O(C_1)|=2$, and then $C_i$ is a circuit for $2\le i\le t$.

Since $|O(C_1)|=2$ and $G$ is connected, we have $|\bar{E}(C_1)|\ge2$: for otherwise, $|V(G)|=|V(C_1)|=2$, Theorem \ref{th:induction step} holds obviously. Hence, by Lemma \ref{lem:x:8'} (2), $L(G)$ has a nontrivial path $P$ with $V(P)=\bar{E}(C_1)$. Now we claim that $|\bar{E}(C_i)|\ge3$ for $2\le i\le t$: if $C_i$ is nontrivial, then we are done; if $C_i$ is an isolated vertex, then by the definition of $EUP_{n}(G)$, $d_G(C_i)\ge3$, our claim holds. By Lemma \ref{lem:x:8} (2), we can find a cycle $C_i'$ in $L(G)$ with $V(C_i')=\bar{E}(C_i)\,(2\le i\le t)$. Let \[ H'=P\cup\big(\bigcup_{i=2}^t C_i'\big).\] We will prove that $H'\in EUP_k(L(G))$.

By Lemma \ref{lem:x:11} and the minimality of $t$, $E(P)\cap E(C_i')=\emptyset$ and $E(C_i')\cap E(C_j')=\emptyset$ for $2\le i,j\le t$ with $i\neq j$, which implies that $|O(H')|=2$. (I)$'$ holds.

Since $P$ is nontrivial, and $V(C_i')=\bar{E}(C_i)\ge3\,(2\le i\le t)$, $H'$ contains no isolated vertex. Note that $\bigcup_{i=3}^{\Delta(G)}V_i(G)\subseteq V(H)$ and $V(H')=\bigcup_{i=1}^t \bar{E}(C_i)$. We have \[\bigcup_{i=3}^{\Delta(L(G))}V_i(L(G))\subseteq V(H').\] Hence, $H'$ satisfies (II).

The details of $H'$ satisfying (III), (IV) and (V)$'$ are almost the same as the proof of Theorem \ref{th:x:induction step} in \cite{X}, so we omit them here. 

It follows that $H'\in EUP_k(L(G))$.

\vspace{.1in}

\textit{Necessity}. Supposing that $EUP_{k}(L(G))\neq\emptyset$. Note that if $EU_{k}(L(G))\neq\emptyset$, then by Theorem \ref{th:x:induction step}, $EU_{k+1}(G)\neq\emptyset$, and hence, $EUP_{k+1}(G)\neq\emptyset$. So we may assume that $EU_{k}(L(G))=\emptyset$, which implies that each subgraph of $L(G)$ in $EUP_{k}(L(G))$ contains exactly two odd vertices.

Let $H$ be a subgraph of $L(G)$ in $EUP_{k}(L(G))$ with a minimum number of isolated vertices. Then $H$ contains no isolated vertices. For otherwise, suppose $C_1=\{e_0\}$ is an isolated vertex of $H$, then by (II), $d_{L(G)}(e_0)\ge3$.  By Lemma \ref{lem:claw-free}, there exist $e_1,e_2\in N_{L(G)}(e_0)$ such that $e_1e_2\in E(L(G))$. Now we construct a subgraph $H_0$ of $L(G)$ as follows.

\[H_0=\left\{\begin{array}{ll}H+\{e_0e_1,e_1e_2,e_2e_0\}&\text{ if }e_1e_2\notin E(H),\\
H+\{e_0e_1,e_2e_0\}-\{e_1e_2\}&\text{ if }e_1e_2\in E(H).
\end{array}\right.\]
Obviously $H_0\in EUP_{k}(L(G))$ has fewer isolated vertices than $H$ has, a contradiction.

Let $H_1,H_2,\dots,H_m$ be the components of $H$, and without loss of generality, let $|O(H_1)|=2$. Since $H$ has no isolated vertices, $H_i$ is an eulerian subgraph of $L(G)$ for $2\le i\le m$. Then $H_1$ can be decomposed into a nontrivial path $P$ and several eulerian subgraphs. Let $P,H_1^1,H_1^2,\dots,H_1^q$ be such a decomposition with $q$ minimized. Then $V(H_1^i)\cap V(H_1^j)=\emptyset$ for $\{i,j\}\subseteq\{1,2,\dots,q\}$ with $i\neq j$.

For the path $P$, by Lemma \ref{lem:x:8'} (1), $G$ has a trail $T'$, such that $L(G[\bar{E}(T')])$ contains $P$, and $P$ contains all elements of $E(T')$. For any eulerian subgraph $H_1^j$ ($1\le j\le q$) or $H_i$ ($i\in\{2,3,\dots,m\}$), by Lemma \ref{lem:x:12}, there exists a subgraph $C_1^j$ or $C_i$ of $G$, respectively, satisfying (1) to (4) of Lemma \ref{lem:x:12}. Let \[C=\big(\bigcup_{i=3}^{\Delta(G)} V_i(G)\big)\cup\big(T\Delta\big(\bigcup_{j=1}^q C_1^j\big)\big)\cup \big(\bigcup_{i=2}^m C_i\big),\] where $T\Delta(\bigcup_{j=1}^q C_1^j)$ is the subgraph of $G$ with vertex set $V(T\cup(\bigcup_{j=1}^q C_1^j))$ and edge set $E(T)\Delta E(\bigcup_{j=1}^q C_1^j)$.

We will prove that $C\in EUP_{k+1}(G)$.

Since $V(H_i)\cap V(H_j)=\emptyset$ for $\{i,j\}\subseteq\{1,2,\dots,m\}$ with $i\neq j$,  $V(H_1^j)\subseteq V(H_1)$ ($1\le j\le q$), and $V(H_1^i)\cap V(H_1^j)=\emptyset$ for $\{i,j\}\subseteq\{1,2,\dots,q\}$ with $i\neq j$, we have $E(C_i)\cap E(C_j)=\emptyset$, $E(C_1^i)\cap E(C_1^j)=\emptyset$, and $E(C_i)\cap E(C_1^j)=\emptyset$. It follows that $d_C(x)\equiv 0\, (\!\!\!\!\mod 2)$ for every $x\in V(C)$ excepting the end-vertices of the path $T$, which implies that $C$ satisfies (I)$'$. Since $C_i$ and $C_1^j$ satisfy Lemma \ref{lem:x:12} (2), $d_G(x)\ge3$ for every $x\in V(C)$ with $d_C(x)=0$. Thus, (II) holds.

Arguing similarly as the proof of Theorem \ref{th:x:induction step}, $C$ satisfies (III), (IV) and (V)$'$. 

It follows that $C\in EUP_{k+1}(G)$.

This completes the proof of Theorem \ref{th:induction step}.\qed

\begin{theorem}\label{th:basis step} \itshape Let $G$ be a connected graph with at least three edges. Then $L^2(G)$ is traceable if and only if $EUP_2(G)\neq\emptyset$.\end{theorem}

\medbreak\noindent \textbf{Proof.} \textit{Sufficiency}. Supposing that $EUP_{2}(G)\neq\emptyset$. Note that if $EU_{2}(G)\neq\emptyset$, then by Theorem \ref{th:x:induction step}, $L^2(G)$ is hamiltonian, we are done. So we may assume that $EU_{2}(G)=\emptyset$, which implies that each subgraph of $G$ in $EUP_{2}(G)$ contains exactly two odd vertices.

We choose an $H\in EUP_{2}(G)$ with a minimum number of components that are denoted by $H_1,H_2,\dots,H_t$, and assume that $|O(H_1)|=2$. Since $H\in EUP_{2}(G)$ and $|E(G)|\ge3$, we have $|\bar{E}(H_1)|\ge2$, and $|\bar{E}(H_i)|\ge3$ for $i\in\{2,3,\dots,t\}$. Then by Lemma \ref{lem:x:8'} (2), we can find a nontrivial path $P$ of $L(G)$ such that $V(P)=\bar{E}(H_1)$. By Lemma \ref{lem:x:8} (2), we can find a cycle $C_i$ of $L(G)$ such that $V(C_i)=\bar{E}(H_i)$, $i\in\{2,3,\dots,t\}$. Let \[T=P\cup\big(\bigcup_{i=2}^t C_i\big).\] By Lemma \ref{lem:x:11} and the minimality of $t$, $P,C_2,C_3,\dots,C_t$ are edge-disjoint. Hence, $T$ is a subgraph of $L(G)$ with exactly 2 odd vertices. Since $d_G(H',H-H')\le 1$ for every subgraph $H'$ of $H$, $T$ is connected.

Note that $H$ satisfies (II), $V(P)=\bar{E}(H_1)$ and $V(C_i)=\bar{E}(H_i)$ for $i\in\{2,3,\dots,t\}$. By the fact that any edge in $G$, which corresponds to a vertex of degree at least 3 in $L(G)$, must incident to a vertex of degree at least 3 in $G$, \[\bigcup_{i=3}^{\Delta(L(G))}V_i(L(G))\subseteq V(T).\]
Since $H\in EUP_{2}(G)$, any branch $b\in B(L(G))$ with $E(b)\cap E(C)=\emptyset$ has length at most 2, and any branch in $B_1(L(G))$ has length at most 1. Then by Lemma \ref{lem:x:10}, $\bar{E}(T)=E(L(G))$, which implies that $T$ is a dominating trail of $L(G)$. Hence, $L^2(G)$ is traceable by Theorem \ref{th:zong}.

\vspace{.1in}

\textit{Necessity}. Supposing that $L^2(G)$ is traceable. By Theorem \ref{th:zong}, $L(G)$ has a dominating trail. Select a dominating trail $T$ of $L(G)$ with a maximum number of vertices of degree at least 3.

\begin{claim}\label{c:2} $\bigcup_{i=3}^{\Delta(L(G))} V_i(L(G))\subseteq V(T)$.\end{claim}

The proof of Claim \ref{c:2} is the same as the proof that $H$ has no isolated vertices in Theorem \ref{th:induction step}, so we omit it here. 

\vspace{.1in}

Then $T$ can be decomposed into a nontrivial path $P$ and several eulerian subgraphs. Let $P,H_1,H_2,\dots,H_q$ be such a decomposition with $q$ minimized. Note that if $T$ is closed, then $T=\cup_{i=1}^q H_i$.

For the path $P$, by Lemma \ref{lem:x:8'} (1), $G$ has a trail $T'$, such that $L(G[\bar{E}(T')])$ contains $P$, and $P$ contains all elements of $E(T')$. For any eulerian subgraph $H_i$ ($1\le i\le q$), by Lemma \ref{lem:x:12}, there exists a subgraph $C_i$ of $G$, satisfying (1) to (4).

Set \[H=\big(\bigcup_{i=3}^{\Delta(G)} V_i(G)\big)\cup\big(T'\Delta\big(\bigcup_{i=1}^q C_i\big)\big),\] where $T'\Delta(\bigcup_{i=1}^q C_i)$ is the subgraph of $G$ with vertex set $V(T'\cup(\bigcup_{i=1}^q C_i))$ and edge set $E(T')\Delta E((\bigcup_{i=1}^q C_i))$. We will prove that $H\in EUP_2(G)$. Before this, we present the following claim.

\begin{claim}\label{c:1} $d_G(x,H)\le1$ for any $x\in \bigcup_{i=3}^{\Delta(G)} V_i(G)$.\end{claim}

\medbreak\noindent \textit{Proof of Claim \ref{c:1}.} If $G$ is either a star or a cycle, then the conclusion holds. For otherwise, then $E_G(x)\cap (\cup_{i=3}^{\Delta(L(G))}V_i(L(G)))\neq\emptyset$ for every vertex $x$ in $\cup_{i=3}^{\Delta(G)}V_i(G)$. Hence, by Claim \ref{c:2}, there exists an edge $e_x$, which is incident to $x$ in $G$, has an endvertex in $H$. Claim \ref{c:1} holds.

\vspace{.1in}

Now we prove $H\in EUP_2(G)$. Obviously, $H$ satisfies (I)$'$. Since $C_i\,(1\le i\le q)$ satisfies (2) of Lemma \ref{lem:x:12}, and by the definition of $H$, (II) holds. By Claim \ref{c:1} and  (3) of Lemma \ref{lem:x:12}, $d_G(H',H-H')\le1$ for every subgraph $H'$ of $H$, thus $H$ satisfies (III). It follows from Lemma \ref{lem:x:10} and $E(L(G))=\bar{E}(T)$ that $|E(b)|\le3$ for $b\in B(G)$ with $E(b)\cap E(H)=\emptyset$, and $|E(b)|\le2$ for $b\in B_1(G)$ with $E(b)\cap E(H)=\emptyset$. $H$ satisfies (IV) and (V)$'$.  Hence, $H\in EUP_2(G)$.

This completes the proof of Theorem \ref{th:basis step}.\qed

Now we prove Theorem \ref{th:the main}.

\medbreak\noindent \textbf{Proof of Theorem \ref{th:the main}.} We proceed by induction on $n$. Theorem \ref{th:basis step} shows that Theorem \ref{th:the main} holds for $n=2$.

Assume, as an inductive hypothesis, that the theorem is true for $n=k>2$, i.e., $L^k(G)$ is traceable if and only if $EUP_k(G)\neq\emptyset$. Now let $n=k+1$. Then $L^{k+1}(G)=L^k (L(G))$ is traceable if and only if $EUP_k(L(G))\neq\emptyset$. Hence, by Theorem \ref{th:induction step},  $EUP_k(L(G))\neq\emptyset$ if and only if $EUP_{k+1}(G)\neq\emptyset$. Theorem \ref{th:the main} holds for $n=k+1$. Thus, the induction succeeds.\qed

Note that Theorem \ref{th:the main} doesn't hold for $n=1$ (see Section \ref{sec:intro}). Hence, when we prove it by induction, the basis step is $n=2$ (Theorem \ref{th:basis step}). 

\section{Applications of Theorem \ref{th:the main}}\label{sec:index}

\hspace*{\parindent} In this section, inspired by the massive upper bounds of $h(G)$, as applications of Theorem \ref{th:the main}, we will present some upper bounds on the hamiltonian path index $h_p(G)$ of a graph $G$. The main idea is to show that $EUP_k(G)\neq\emptyset$, and then by Theorem \ref{th:the main}, $h_p(G)\le k$. 

Note that the hamiltonian index $h(G)$ exists for any connected graph $G$ other than a path and $h_p(G)\le h(G)$. The hamiltonian path index $h_p(G)$ exists for any connected graph.

Comparing the definition of hamiltonian with traceable, we know that being hamiltonian is stronger than being traceable. Then one may believe that the former needs less iterated steps, and hence, $h_p(G)< h(G)$. Unfortunately, this is not true by the fact that $h_p(K_{1,n-1})=h(K_{1,n-1})=1$ ($n\ge3$). Moreover, Fig. 2 shows a graph $G$ with $h_p(G)=h(G)=k$: one can check that the unique cycle of $G$ is an element in $EU_k(G)$, but $EUP_{k-1}(G)=\emptyset$ by the fact that any element in  $EUP_{k-1}(G)$ cann't contain all the three pendent paths with length $k$. Hence, our trivial bound $h_p(G)\le h(G)$ is best possible.

\vspace{.2in}
\begin{center}
\begin{tikzpicture}[line width=.6pt]
\draw (3.1,1) arc (0:360:1.1);
\filldraw (2.955,1.55) circle [radius=1.8pt]
            (2,2.1) circle [radius=1.8pt]
            (1.045,1.55) circle [radius=1.8pt]
            (1.045,0.45) circle [radius=1.8pt]
            (2,-0.1) circle [radius=1.8pt]
            (2.955,0.45) circle [radius=1.8pt];
\filldraw (2,2.7) circle [radius=1.8pt]
            (2,3.1) circle [radius=.8pt]
            (2,3.3) circle [radius=.8pt]
            (2,3.5) circle [radius=.8pt]
            (2,4) circle [radius=1.8pt]
            (2,4.6) circle [radius=1.8pt];
\draw (2,2.1) -- (2,2.9);
\draw (2,3.8) -- (2,4.6); 
\draw [gray,decorate,decoration={brace,mirror,amplitude=0.15cm}] (2.2,2.15) - - (2.2,4.55);
\draw[white] (2.8,2.1) -- (2.8,4.6) node[gray,midway,above,sloped] {\tiny a path of};
\draw[white] (3.05,2.1) -- (3.05,4.6) node[gray,midway,above,sloped] {\tiny length $k$};           
\filldraw (0.445,0.45) circle [radius=1.8pt]
            (0.045,0.45) circle [radius=.8pt]
            (-.155,0.45) circle [radius=.8pt]
            (-.355,0.45) circle [radius=.8pt]
            (-.755,0.45) circle [radius=1.8pt]
            (-1.355,0.45) circle [radius=1.8pt]; 
\draw (-1.355,0.45) -- (-.555,0.45);
\draw (.245,0.45) -- (1.045,0.45);
\draw [gray,decorate,decoration={brace,mirror,amplitude=0.15cm}] (-1.35,.25) - - (1.04,.25);  
\draw[gray] (-.155,.1) node[below]{\tiny a path of length $k$};                     
\filldraw (3.555,0.45) circle [radius=1.8pt]
            (3.955,0.45) circle [radius=.8pt]
            (4.155,0.45) circle [radius=.8pt]
            (4.355,0.45) circle [radius=.8pt]
            (4.755,0.45) circle [radius=1.8pt]
            (5.355,0.45) circle [radius=1.8pt];            
\draw (2.955,0.45) -- (3.755,0.45);
\draw (4.555,0.45) -- (5.355,0.45);
\draw [gray,decorate,decoration={brace,mirror,amplitude=0.15cm}] (3.005,.25) - - (5.305,.25);  
\draw[gray] (4.155,.1) node[below]{\tiny a path of length $k$};
\end{tikzpicture}

\vspace{.15in}
\small Fig. 2. A graph $G$ with $h_p(G)=h(G)=k$.
\end{center}

For a graph $G$, let $MT^*(G)$ be a trail of $G$ with the most number of vertices, and in this sence, with the least number of vertices in $\bigcup_{i=3}^{\Delta(G)}V_i(G)$. Denote $mt^*(G)=|V(MT^*(G))|$ and $d_{\ge3}^*(G)=|\bigcup_{i=3}^{\Delta(G)}V_i(G)\backslash V(MT^*(G))|$.

\begin{theorem}\label{th:b1} \itshape Let $G$ be a connected graph of order $n$. Then $h_p(G)\le n-mt^*(G)-d_{\ge3}^*(G)+2$.\end{theorem}

\medbreak\noindent \textbf{Proof.} Since $G$ is connected, Theorem \ref{th:b1} holds for $|E(G)|<3$ trivially. So we may assume that $|E(G)|\ge 3$.

Let $MT^*(G)$ be a trail of $G$ satisfying the hypotheses above. Denote $k=n-mt^*(G)-d_{\ge3}^*(G)+2$. Note that $k\ge2$. By Theorem \ref{th:the main}, it suffices to prove that $EUP_k(G)\neq\emptyset$.

Let $$H=MT^*(G)\cup\big(\bigcup_{i=3}^{\Delta (G)}V_i(G)\big).$$ We will prove that $H\in EUP_k(G)$.

By the definition of $H$, $H$ satisfies (I)$'$ and (II). Note that $|V(G)|-|V(H)|=n-(mt^*(G)+d_{\ge3}^*(G))=k-2$. Then $d_G(H_1,H-H_1)\le k-1$ for every subgraph $H_1$ of $H$, $|E(b)|\le k-1\,(<k+1)$ for every branch $b\in B(G)$ with $E(b)\cap E(H)=\emptyset$, and $|E(b)|\le k-2\,(<k)$ for every branch $b\in B_1(G)$  with $E(b)\cap E(H)=\emptyset$. Hence, $H$ satisfies (III), (IV) and (V)$'$. Theorem \ref{th:b1} holds.
\qed

The bound of $h_p(G)$ in Theorem \ref{th:b1} is sharp. Fig. 3 shows a graph $G$ with $h_p(G)=n-mt^*(G)-d_{\ge3}^*(G)+2$, where $s$, $t$ are positive integers and $t\ge s+5$ in the figure. Since $t\ge s+5$, $MT^*(G)$ is the path $x_1x_2\cdots x_tw_1y_t\cdots y_2y_1$, and then $mt^*(G)=2t+1$ and $d_{\ge3}^*(G)=4$. Note that $n=|V(G)|=2t+s+5$. On the one hand, by Theorem \ref{th:b1}, we have $h_p(G)\le n-mt^*(G)-d_{\ge3}^*(G)+2=s+2$. On the other hand, we will explain $h_p(G)\ge s+2$ in the following. For $k\le s+1$, the $k$-th iterated line graph $L^k(G)$ is also illustrated in Fig. 3, where the gray ellipse and triangle are the nontrivial hamiltonian subgraph $S_1$ and $S_2$ of $L^k(G)$, respectively. Note that $k\le s+1$. We have $s-k+2\ge1$, which means that either $|V(S_1)\cap V(S_2)|=1$ (when $s-k+2=1$), or $d_{L^k(G)}(S_1,S_2)\ge1$ (when $s-k+2\ge2$). In both cases, $L^k(G)$ is not traceable. Hence, $h_p(G)\ge s+2$.

\vspace{.1in}
\begin{center}
\begin{tikzpicture}[line width=.6pt]
\draw (0,0) -- (1.1,0);
\draw (1.9,0) -- (4.1,0);
\draw (4.9,0) -- (6,0);
\draw (3,0) -- (3,.8);
\draw (3,1.6) -- (3,3.8);
\filldraw (0,0) circle [radius=1.8pt]
            (.8,0) circle [radius=1.8pt]
            (1.3,0) circle [radius=.8pt]
            (1.5,0) circle [radius=.8pt]
            (1.7,0) circle [radius=.8pt]
            (2.2,0) circle [radius=1.8pt]
            (3,0) circle [radius=1.8pt]
            (3.8,0) circle [radius=1.8pt]
            (4.3,0) circle [radius=.8pt]
            (4.5,0) circle [radius=.8pt]
            (4.7,0) circle [radius=.8pt]
            (5.2,0) circle [radius=1.8pt]
            (6,0) circle [radius=1.8pt];
\filldraw (3,.6) circle [radius=1.8pt]
            (3,1) circle [radius=.8pt]
            (3,1.2) circle [radius=.8pt]
            (3,1.4) circle [radius=.8pt]
            (3,1.8) circle [radius=1.8pt]
            (3,2.4) circle [radius=1.8pt];
\draw (3,3.1) circle [radius=.7];
\filldraw (3,3.8) circle [radius=1.8pt]
            (2.3,3.1) circle [radius=1.8pt]
            (3.7,3.1) circle [radius=1.8pt];
\draw (2.3,3.1) -- (3.7,3.1);
\draw [gray,decorate,decoration={brace,mirror,amplitude=0.15cm}] (0.05,-.2) - - (2.15,-.2);
\draw [gray,decorate,decoration={brace,mirror,amplitude=0.15cm}] (3.85,-.2) - - (5.95,-.2);
\draw [gray,decorate,decoration={brace,mirror,amplitude=0.15cm}] (3.2,0.65) - - (3.2,1.75);
\draw[gray] (1.1,-.35) node[below]{\tiny a path of $t$ vertices};
\draw[gray] (4.9,-.35) node[below]{\tiny a path of $t$ vertices};
\draw[white] (3.8,.6) -- (3.8,1.7) node[gray,midway,above,sloped] {\tiny a path of};
\draw[white] (4,.6) -- (4,1.7) node[gray,midway,above,sloped] {\tiny $s$ vertices};
\draw (0,0) node[above]{$x_1$};
\draw (.8,0) node[above]{$x_2$};
\draw (2.2,0) node[above]{$x_t$};
\draw (6,0) node[above]{$y_1$};
\draw (5.2,0) node[above]{$y_2$};
\draw (3.8,0) node[above]{$y_t$};
\draw (3,.6) node[left]{$z_1$};
\draw (3,1.8) node[left]{$z_s$};
\draw (3,0) node[below]{$w_1$};
\draw (3,2.2) node[right]{$w_2$};
\draw (2.3,3.1) node[left]{$w_3$};
\draw (3,3.75) node[above]{$w_4$};
\draw (3.7,3.1) node[right]{$w_5$};
\draw (3,-.85) node[below]{$G$};

\draw [fill=gray!20,rounded corners=0.3cm](9.4,-.14) - - (10.5,1.2)  - - (11.6,-.14)  - - cycle;
\filldraw[fill=gray!20] (10.5,3) ellipse (.5 and .87);
\draw (7.5,0) -- (8.6,0);
\draw (12.4,0) -- (13.5,0);
\draw (9.4,0) -- (9.7,0);
\draw (11.3,0) -- (11.6,0);
\draw (10.5,.9) -- (10.5,1.2);
\draw (10.5,2) -- (10.5,2.3);
\filldraw (7.5,0) circle [radius=1.8pt]
            (8.3,0) circle [radius=1.8pt]
            (8.8,0) circle [radius=.8pt]
            (9,0) circle [radius=.8pt]
            (9.2,0) circle [radius=.8pt]
            (9.7,0) circle [radius=1.8pt]
            (11.3,0) circle [radius=1.8pt]
            (11.8,0) circle [radius=.8pt]
            (12,0) circle [radius=.8pt]
            (12.2,0) circle [radius=.8pt]
            (12.7,0) circle [radius=1.8pt]
            (13.5,0) circle [radius=1.8pt];
\filldraw (10.5,.9) circle [radius=1.8pt]
            (10.5,1.4) circle [radius=.8pt]
            (10.5,1.6) circle [radius=.8pt]
            (10.5,1.8) circle [radius=.8pt]
            (10.5,2.3) circle [radius=1.8pt];
\draw [gray,decorate,decoration={brace,mirror,amplitude=0.15cm}] (7.55,-.2) - - (9.65,-.2);
\draw [gray,decorate,decoration={brace,mirror,amplitude=0.15cm}] (11.35,-.2) - - (13.45,-.2);
\draw [gray,decorate,decoration={brace,mirror,amplitude=0.15cm}] (10.6,.95) - - (10.6,2.25);
\draw[gray] (8.6,-.35) node[below]{\tiny a path of $t-k+1$ vertices};
\draw[gray] (12.4,-.35) node[below]{\tiny a path of $t-k+1$ vertices};
\draw[white] (11.2,.95) -- (11.2,2.25) node[gray,midway,above,sloped] {\tiny a path of};
\draw[white] (11.45,.95) -- (11.45,2.25) node[gray,midway,above,sloped] {\tiny $s-k+2$ vertices};
\draw (10.5,3) node{$S_1$};
\draw (10.5,.35) node{$S_2$};
\draw (10.5,-.85) node[below]{$L^k(G)$ with $k\le s+1$};
\end{tikzpicture}

\vspace{.1in}
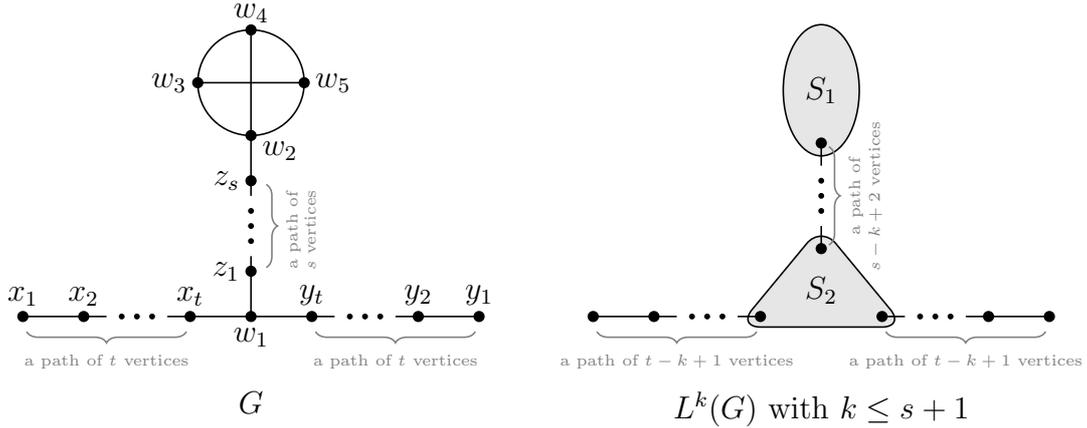
\small Fig. 3. A graph $G$ with $h_p(G)=n-mt^*(G)-d_{\ge3}^*(G)+2$ 

\quad\quad$\,\,$ and its iterated line graph $L^k(G)$  with $k\le s+1$.
\end{center}

By Theorem \ref{th:b1}, we can obtain the following corollary.

\begin{corollary}\label{cor1:b1} \itshape Let $G$ be a connected graph of order $n$. Then $h_p(G)\le \max\{1,n-mt^*(G)\}$.\end{corollary}

\medbreak\noindent \textbf{Proof.} Let $MT^*(G)$ be a maximum trail of $G$, and in this sence, with the least number of vertices in $\bigcup_{i=3}^{\Delta(G)}V_i(G)$, and let $k=n-mt^*(G)$. If $k\le1$, then $MT^*(G)$ is a dominating trail of $G$. By Theorem \ref{th:zong}, $L(G)$ is traceable. Hence, $h_p(G)\le 1$. So we may assume that $k\ge2$. Now the proof is divided into three cases.

\medbreak\noindent \textbf{Case 1.} $d_{\ge3}^*(G)=0$. 

\medbreak Let $H=MT^*(G)$. We will prove that $H\in EUP_k(G)$, and then $h_p(G)\le k$. Obviously, $H$ satisfies (I)$'$, (II) and (III). Note that $d_{\ge3}^*(G)=0$. $G-H$ has exactly $k$ vertices. Then $|E(b)|\le k+1$ for every branch $b\in B(G)$ with $E(b)\cap E(H)=\emptyset$, and $|E(b)|\le k$ for every branch $b\in B_1(G)$  with $E(b)\cap E(H)=\emptyset$. Hence, $H$ satisfies (IV) and (V)$'$. 

\medbreak\noindent \textbf{Case 2.} $d_{\ge3}^*(G)=1$.

Let $v$ be the vertex of degree at least 3 in $\bigcup_{i=3}^{\Delta(G)}V_i(G)\backslash V(MT^*(G))$. If $|N_G(v)|=1$, then by $d_{\ge3}^*(G)=1$, the neighbour of $v$ belongs to $MT^*(G)$. Hence, the union of $MT^*(G)+v$ and two multiple edges incident to $v$ is a longer trail than $MT^*(G)$, contrary to the maximality of $MT^*(G)$. Then we may assume that $|N_G(v)|\ge2$.

Let $H=MT^*(G)+v$. We will prove that $H\in EUP_k(G)$. Obviously, $H$ satisfies (I)$'$ and (II). Arguing similarly as the proof of Case 1,  $H$ satisfies (IV) and (V)$'$. It remains to prove that $d_G(MT^*(G),v)\le k-1$. This holds  by the fact that $G-H$ has exactly $k-1$ vertices, $|N_G(v)|\ge2$, and the shortest path between $v$ and $MT^*(G)$ contains only one neighbour of $v$.

\medbreak\noindent \textbf{Case 3.} $d_{\ge3}^*(G)\ge2$.

\medbreak By Theorem \ref{th:b1}, $h_p(G)\le n-mt^*(G)-d_{\ge3}^*(G)+2$. Then by $d_{\ge3}^*(G)\ge2$, we have $h_p(G)\le n-mt^*(G)$. 

This completes the proof of Corollary \ref{cor1:b1}. \qed

By the sharpness of $h_p(G)\le n-mt^*(G)-d_{\ge3}^*(G)+2$ and the proof of Corollary \ref{cor1:b1}, we know that the bound $h_p(G)\le \max\{1,n-mt^*(G)\}$ is sharp when $d_{\ge3}^*(G)\le2$.

The \textit{diameter} of a graph $G$, denoted by $\mathrm{diam}(G)$, is the greatest distance between two vertices of G. Note that $\mathrm{diam} (G)+1\le mt^*(G)$. The following corollary is obvious.

\begin{corollary}\label{cor2:b1} \itshape Let $G$ be a connected graph of order $n$. Then $h_p(G)\le \max\{1,n-\mathrm{diam}(G)-1\}$.\end{corollary}

For a graph $G$, let $d'_G(v)=|N_G(v)|$ and $\Delta'(G)=\max\{d'_G(v)\,:\,v\in V(G)\}$. Note that if $G$ is simple, then $d'_G(v)=d_G(v)$ and $\Delta'(G)=\Delta(G)$. Let $d_{\ge3}^{**}(G)=\max\{|(\bigcup_{i=3}^{\Delta (G)}V_i(G))\backslash N_G(v)|\,:\,v\in V(G) \text{ and } |N_G(v)|=\Delta'(G)\}$. A cycle $C$ of $G$ is called \textit{pendent} if $|V(C)\cap (\bigcup_{i=3}^{\Delta (G)}V_i(G))|=1$. See Fig. 4 (a) for illustrations of pendent cycles. Let $PC(G)$ be the set of pendent cycles of $G$. 

\begin{theorem}\label{th:b2} \itshape Let $G$ be a connected graph of order $n$. Then $$h_p(G)\le \bigg\lfloor\frac{\,n-\Delta'(G)-d_{\ge3}^{**}(G)\,}{\,3\,}\bigg\rfloor+3.$$\end{theorem}

\medbreak\noindent \textbf{Proof.} Let $k=\lfloor (n-\Delta'(G)-d_{\ge3}^{**}(G))/3\rfloor+3$, and $v$ a vertex of $G$ with $|N_G(v)|=\Delta'(G)$ and $|(\bigcup_{i=3}^{\Delta (G)}V_i(G))\backslash N_G(v)|$ maximized. Note that $k\ge3$. By Theorem \ref{th:the main}, it suffices to prove that $EUP_k(G)\neq\emptyset$.

If $\Delta'(G)\le2$, then $G$ is traceable. Hence, $h_p(G)=0$, and the bound holds trivially. Now we may assume that $\Delta'(G)\ge3$. Then $B(G)$, the set of branches of $G$, has at least two elements.

Let $b_1\in B(G)$ be a branch with $|V(b_1)\cap (V_1(G)\cup V_2(G))|$ maximized, and $b_2\in B(G)\backslash \{b_1\}$ a branch with $|V(b_2)\cap (V_1(G)\cup V_2(G))|$ maximized. Since $G$ is connected, we can find a trail $T_1$ (may be trivial) which connects $b_1$ and $b_2$. By the fact that the internal vertices (if any) of $b_1$ and $b_2$ have degree 2 in $G$, $T=T_1\cup b_1\cup b_2$ is a trail of $G$.

Let \[H=\big(\bigcup_{i\ge3}^{\Delta(G)}V_i(G)\big)\cup T\cup PC(G).\] We will prove that $H\in EUP_k(G)$. 

Obviously, $H$ satisfies (I)$'$ and (II). By the choice of $b_1$ and $b_2$, each of the other branches of $G$ has at most $\lfloor (n-\Delta'(G)-d_{\ge3}^{**}(G))/3\rfloor+1$ vertices in $V_1(G)\cup V_2(G)$, where the $+1$ is necessary since each branch may have at most one neighbour of $v$. Then $|E(b)|\le \lfloor (n-\Delta'(G)-d_{\ge3}^{**}(G))/3\rfloor+2=k-1$ for every branch $b\in B(G)\backslash \{b_1,b_2\}$, and hence, (III), (IV) and (V)$'$ hold.

This completes the proof of Theorem \ref{th:b2}.\qed

The bound of $h_p(G)$ in Theorem \ref{th:b2} is sharp. Fig. 4 (b) shows a graph $G$ with $n=3s+16$, $\Delta'(G)=6$ and $d_{\ge3}^{**}(G)=13$, where $s$ is a positive integer and the three gray cycles are induced $K_4$. Then by Theorem \ref{th:b2}, $h_p(G)\le  \lfloor (n-\Delta'(G)-d_{\ge3}^{**}(G))/3\rfloor+3=s+2$. Note that $G$ has 3 branches of length $s+1$. $L^k(G)$ is not traceable when $k<s+2$. Then $h_p(G)\ge s+2$, and hence, $h_p(G)=s+2$, which implies the sharpness of the upper bound in Theorem \ref{th:b2}.

If $G$ is a connected simple graph of order $n$, then Theorem \ref{th:b2} implies that $h_p(G)\le \lfloor (n-\Delta(G)-d_{\ge3}^{**}(G))/3\rfloor+3$.

\vspace{.2in}
\begin{center}
\begin{tikzpicture}[line width=.6pt]
\draw (0,0) arc (170:110:4);
\filldraw (0,0) circle [radius=1.8pt]
            (.182,.675) circle [radius=1.8pt]
            (.475,1.31) circle [radius=1.8pt]
            (.874,1.88) circle [radius=1.8pt]
            (1.37,2.375) circle [radius=1.8pt]
            (1.945,2.772) circle [radius=1.8pt]
            (2.57,3.065) circle [radius=1.8pt];
\draw (1.37,2.375) arc (-50:310:.6);
\filldraw (1.445,3.22) circle [radius=1.8pt]
            (.6,3.295) circle [radius=1.8pt]
            (.53,2.45) circle [radius=1.8pt];
\filldraw (1.6,.75) circle [radius=1.8pt];
\draw (.475,1.31) to[bend left=45] (1.6,.75);
\draw (.475,1.31) to[bend left=-45] (1.6,.75);
\draw (1.285,-.55) node[below]{(a)};

\draw (7,1.2) circle [radius=1pt];
\draw[fill=gray!20] (7,3.6) circle [radius=.6];
\draw[fill=gray!20] (4.9215,0) circle [radius=.6];
\draw[fill=gray!20] (9.0785,0) circle [radius=.6];
\filldraw (7,1.2) circle [radius=1.8pt];
\filldraw (7,.6) circle [radius=1.8pt];
\filldraw (7,1.8) circle [radius=1.8pt];

\filldraw (7.5196,.9) circle [radius=1.8pt];
\filldraw (7.5196,1.5) circle [radius=1.8pt];
\filldraw (6.4804,.9) circle [radius=1.8pt];
\filldraw (6.4804,1.5) circle [radius=1.8pt];
\draw (7,1.2) -- (7,.6);
\draw (7,1.2) -- (7.5196,1.5);
\draw (7,1.2) -- (6.4804,1.5);
\draw (7,1.2) -- (7,1.95);
\draw (7,3.2) -- (7,2.45);

\filldraw (7,3.2) circle [radius=1.8pt];
\filldraw (7,2.6) circle [radius=1.8pt];

\filldraw (7,2.05) circle [radius=.8pt];
\filldraw (7,2.2) circle [radius=.8pt];
\filldraw (7,2.35) circle [radius=.8pt];

\filldraw (5.2679,.2) circle [radius=1.8pt];
\filldraw (5.7875,.5) circle [radius=1.8pt];
\filldraw (6.1339,.7) circle [radius=.8pt];
\filldraw (6.2638,.775) circle [radius=.8pt];
\filldraw (6.004,.625) circle [radius=.8pt];

\filldraw (8.7321,.2) circle [radius=1.8pt];
\filldraw (8.2125,.5) circle [radius=1.8pt];
\filldraw (7.8661,.7) circle [radius=.8pt];
\filldraw (7.7362,.775) circle [radius=.8pt];
\filldraw (7.996,.625) circle [radius=.8pt];

\draw (7,1.2) -- (6.3505,.825);
\draw (7,1.2) -- (7.6495,.825);
\draw (5.2679,.2) -- (5.9174,.575);
\draw (8.7321,.2) -- (8.0826,.575);
\draw (4.9215,0) node{$K_4$};
\draw (9.0785,0) node{$K_4$};
\draw (7,3.6) node{$K_4$};

\draw[white] (5.6575,.7) -- (6.3504,1.1) node[gray,midway,above,sloped] {\tiny $s$ vertices};
\draw[white] (8.3425,.7) -- (7.6496,1.1) node[gray,midway,above,sloped] {\tiny $s$ vertices};
\draw [gray,decorate,decoration={brace,mirror,amplitude=0.15cm}] (7.2,1.8) - - (7.2,2.6);
\draw [gray,decorate,decoration={brace,amplitude=0.15cm}] (5.7075,.65) - - (6.4004,1.05);
\draw [gray,decorate,decoration={brace,mirror,amplitude=0.15cm}] (8.2925,.65) - - (7.5996,1.05);
\draw (8,2.2) node[gray]{\tiny $s$ vertices};
\draw (7,-.55) node[below]{(b)};
\end{tikzpicture}

\vspace{.1in}
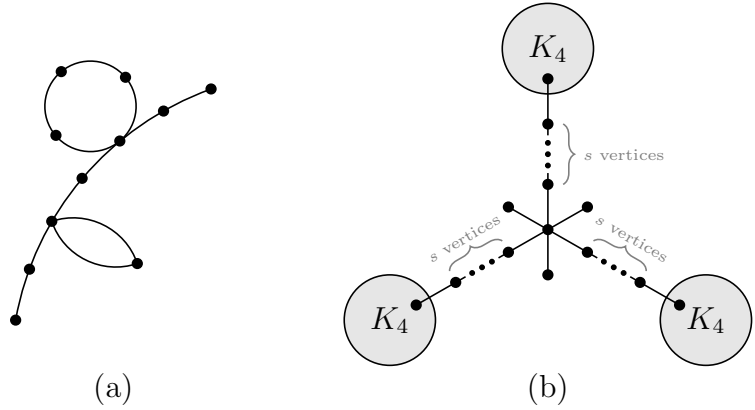
\small Fig. 4. (a) A graph with 2 pendent cycles; and (b) A graph $G$ with  

$h_p(G)= \lfloor (n-\Delta'(G)-d_{\ge3}^{**}(G))/3\rfloor+3$.
\end{center} 

\medbreak In \cite{Sa}, Sara\u{z}in proved the following upper bound of $h(G)$. 

\begin{theorem}(Sara\u{z}in, \cite{Sa})\label{th:sa} \itshape Let $G$ be a connected simple graph of order $n$. If $\Delta(G)\ge3$, then $h(G)\le n-\Delta(G)$.\end{theorem} 

The following corollary, obtained by Theorem \ref{th:b2} immediately, implies that when considering $h_p(G)$, the upper bound in Theorem \ref{th:sa} can be improved evidently.

\begin{corollary}\label{cor:b2} \itshape Let $G$ be a connected simple graph of order $n$. Then $$h_p(G)\le \bigg\lfloor\frac{\,n-\Delta(G)\,}{\,3\,}\bigg\rfloor+3.$$\end{corollary}

%

\end{document}